\documentclass[11pt]{article}
\usepackage{amsmath, amssymb, amsthm, verbatim,enumerate,bbm, mathtools}
\usepackage{indentfirst}

\date{\today}
\allowdisplaybreaks

\theoremstyle{plain}

\newtheorem{theorem}{Theorem}

\newtheorem{corollary}[theorem]{Corollary}

\title{Lower bounds for multicolor Ramsey numbers}
\author{David Conlon\thanks{Department of Mathematics, California Institute of Technology, CA 91125, USA. Email:  {\tt dconlon@caltech.edu}.} \and Asaf Ferber\thanks{Department of Mathematics, University of California, Irvine, CA 92697, USA.
Email: {\tt asaff@uci.edu}.
Research supported in part by NSF grants DMS-1954395 and DMS-1953799.}}
\date{}

\begin{document}

\maketitle

\begin{abstract}
We give an exponential improvement to the lower bound on diagonal Ramsey numbers for any fixed number of colors greater than two.
\end{abstract}

\section{Introduction}

The Ramsey number $r(t; \ell)$ is the smallest natural number $n$ such that every $\ell$-coloring of the edges of the complete graph $K_n$ contains a monochromatic $K_t$. For $\ell = 2$, the problem of determining $r(t) := r(t;2)$ is arguably one of the most famous in combinatorics. The bounds
\[\sqrt{2}^t < r(t) < 4^t\]
have been known since the 1940s, but, despite considerable interest, only lower-order improvements \cite{Con09, Sah20, Spe75} have been made to either bound. In particular, the lower bound $r(t) > (1 + o(1))\frac{t}{\sqrt{2} e} \sqrt{2}^t$, proved by Erd\H{o}s~\cite{Erd47} as one of the earliest applications of the probabilistic method, has only been improved~\cite{Spe75} by a factor of $2$ in the intervening 70 years. 

If we ignore lower-order terms, the best known upper bound for $\ell \geq 3$ is $r(t; \ell) < \ell^{\ell t}$, proved through a simple modification of the Erd\H{o}s--Szekeres neighborhood-chasing argument~\cite{ESz35} that yields $r(t) < 4^t$. For $\ell = 3$, the best lower bound, $r(t; 3) > \sqrt{3}^{t}$, again comes from the probabilistic method. For higher $\ell$, the best lower bounds come from the simple observation of Lefmann~\cite{L87} that
\[r(t; \ell_1 + \ell_2) - 1 \geq (r(t; \ell_1) - 1)(r(t;\ell_2) - 1).\]
To see this, we blow up an $\ell_1$-coloring of $K_{r(t;\ell_1) - 1}$ with no monochromatic $K_t$ so that each vertex set has order $r(t; \ell_2) - 1$ and then color each of these copies of $K_{r(t;\ell_2) - 1}$ separately with the remaining $\ell_2$ colors so that there is again no monochromatic $K_t$.
By using the bounds $r(t; 2) - 1 \geq 2^{t/2}$ and $r(t;3) - 1 \geq 3^{t/2}$, we can repeatedly apply this observation to conclude that
\[r(t; 3k) > 3^{kt/2}, \qquad r(t; 3k+1) > 2^t 3^{(k-1)t/2}, \qquad r(t; 3k +2) > 2^{t/2} 3^{kt/2}.\]
Our main result is an exponential improvement to all these lower bounds for three or more colors.

Our principal contribution is the following theorem, proved via a construction which is partly deterministic and partly random. The deterministic part shares some characteristics with a construction of Alon and Krivelevich~\cite{AK97}, in that we consider a graph whose vertices are vectors over a finite field where adjacency is determined by the value of their scalar product, while randomness comes in through both random coloring and random sampling.

\begin{theorem}
\label{technical}
For any prime $q$, $r(t; q+1) > 2^{t/2} q^{3t/8 + o(t)}$.
\end{theorem}

In particular, the cases $q=2$ and $q= 3$ yield exponential improvements over the previous bounds for $r(t; 3)$ and $r(t;4)$, both of which came from the probabilistic method (in fact, Lefmann's observation gives an additional polynomial factor in the four-color case, but this is of lower order than the exponential improvements that are our concern). 

\begin{corollary}
$r(t;3) > 2^{7t/8 + o(t)}$ and $r(t;4) > 2^{t/2} 3^{3t/8 + o(t)}$.
\end{corollary}

For the sake of comparison, we note that the improvement for three colors is from $1.732^t$ to $1.834^t$, while, for four colors, it is from $2^t$ to $2.135^t$. Improvements for all $\ell \geq 5$ now follow from repeated applications of Lefmann's observation, yielding
\[r(t; 3k) > 2^{7kt/8 + o(t)}, \qquad r(t; 3k+1) > 2^{7(k-1)t/8 + t/2} 3^{3t/8 + o(t)}, \qquad r(t; 3k +2) > 2^{7kt/8 + t/2 + o(t)},\]
where we used, for instance,
\[r(t; 3k+1) - 1 \geq (r(t; 3(k-1)) - 1)( r(t;4) - 1) \geq (r(t; 3) - 1)^{k-1}( r(t;4) - 1).\]

\section{Proof of Theorem \ref{technical}} 

Let $q$ be a prime. Suppose $t \neq 0 \bmod{q}$ and let $V\subseteq \mathbb{F}_q^t$ be the set consisting of all vectors $v\in \mathbb{F}_q^t$ for which $\sum_{i=1}^t v^2_i=0 \bmod{q}$, noting that $q^{t-2} \leq |V|\leq q^t$. Here the lower bound follows from observing that we may pick $v_1, \dots, v_{t-2}$ arbitrarily and, since every element in $\mathbb{F}_q$ can be written as the sum of two squares, there must then exist at least one choice of $v_{t-1}$ and $v_t$ such that $v_{t-1}^2 + v_t^2 = - \sum_{i=1}^{t-2} v_i^2$.

We will first color all the pairs $\binom{V}{2}$ and then define a coloring of $E(K_n)$ by restricting our attention to a random sample of $n$ vertices in $V$. Formally: 

\paragraph{Coloring all pairs in $\binom{V}{2}$.} 

 For every pair $uv\in \binom{V}{2}$, we define its color $\chi(uv)$ according to the following rules: 
\begin{itemize}
    \item If $u\cdot v=i \bmod{q}$ and $i\neq 0$, then set $\chi(uv)=i$. 
    \item Otherwise, choose $\chi(uv)\in \{q,q+1\}$ uniformly at random, independently of all other pairs.
\end{itemize}

\paragraph{Mapping $[n]$ into $V$.}
Take a random injective map $f:[n]\rightarrow V$ and define the color of every edge $ij$ as $\chi(f(i)f(j))$.

\vspace{4mm}

Our goal is to upper bound the orders of the cliques in each color class. 

\paragraph{Colors $1\leq i\leq q-1$.} There are no $i$-monochromatic cliques of order larger than $t$ for any $1\leq i \leq q-1$. Indeed, suppose that $v_1,\ldots,v_s$ form an $i$-monochromatic clique. We will try to show that they are linearly independent and, therefore, that there are at most $t$ of them. 
To this end, suppose that 
$$u:=\sum_{j=1}^s \alpha_j v_j=\bar{0}$$ 
and we wish to show that $\alpha_j=0\bmod{q}$ for all $j$. 
 Observe that since $v_j\cdot v_j=0\bmod{q}$ for all $j$ (our ground set $V$ consists only of such vectors) and $v_k\cdot v_{j}=i\bmod{q}$ for each $k\neq j$, by considering all the products $u\cdot v_j$, we obtain that the vector $\bar{\alpha}=(\alpha_1,\ldots,\alpha_s)$ is a solution to 
 $$M\bar{\alpha}=\bar{0}$$ 
 with $M=iJ-iI$, where $J$ is the $s \times s$ all $1$ matrix and $I$ is the $s \times s$ identity matrix. In particular, we obtain that the eigenvalues of $M$ (over $\mathbb{Z}$) are $is-i$ with multiplicity $1$ and $-i$ with multiplicity $s-1$. Therefore, if $s \neq 1 \bmod{q}$, the matrix is also non-singular over $\mathbb{Z}_q$, implying that $\bar{\alpha}=0$, as required. On the other hand, if $s = 1 \bmod{q}$, we can apply the same argument with $v_1, \dots, v_{s-1}$ to conclude that $s-1 \leq t$. But, we cannot have $s - 1 = t$, since this would imply that $t = 0 \bmod{q}$, contradicting our assumption. Therefore, we may also conclude that $s \leq t$ in this case.

\paragraph{Colors $q$ and $q+1$.} We call a subset $X\subseteq V$ a \emph{potential clique} if $|X|=t$ and $u\cdot v=0\bmod{q}$ for all $u,v\in X$. Given a potential clique $X$, we let $M_X$ be the $t \times t$ matrix whose rows consist of all the vectors in $X$. Observe that $M_X\cdot M_X^T=0$, where we use the fact that each vector is self-orthogonal. First we wish to count the number of potential cliques and later we will calculate the expected number of cliques that survive after we color randomly and restrict to a random subset of order $n$. 

Suppose that $X$ is a potential clique and let $r:=\textrm{rank}(X)$ be the rank of the vectors in this clique, noting that $r \leq t/2$, since the dimension of any isotropic subspace of $\mathbb{F}_q^t$ is at most $t/2$. By assuming that the first $r$ elements are linearly independent, the number of ways to build a potential clique $X$ of rank $r$ is upper bounded by
\[
   \left(\prod_{i=0}^{r-1}q^{t-i}\right) \cdot q^{(t-r)r} =q^{tr - \binom{r}{2} + tr - r^2}
   =q^{2tr - \frac{3r^2}{2} +\frac{r}{2}}.\]
Indeed, suppose that we have already chosen the vectors $v_1,\ldots,v_s \in X$ for some $s<r$. Then, letting $M_s$ be the $s\times t$ matrix with the $v_i$ as its rows, we need to choose $v_{s+1}$ such that $M_s\cdot v_{s+1}=\bar{0}$. Since the rank of $M_s$ is assumed to be $s$, there are exactly $q^{t-s}$ choices for $v_{s+1}$ in $\mathbb{F}_q^t$ and, therefore, at most that many choices for $v_{s+1}\in V$. If, instead, $s\geq r$, then we need to choose a vector $v_{s+1}\in \textrm{span}\{v_1,\ldots,v_r\}$ and there are at most $q^r$ such choices in $V$. 

Now observe that the function $2tr - \frac{3r^2}{2} +\frac{r}{2}$ appearing in the exponent of the expression above is increasing up to $r = \frac{2t}{3} + \frac{1}{6}$, so the maximum occurs at $t/2$. Therefore, by plugging this into our estimate and summing over all possible ranks, we see that the number $N_t$ of potential cliques in $V$ is upper bounded by $q^{\frac{5 t^2}{8} + o(t^2)}$.

The probability that a potential clique becomes monochromatic after the random coloring is $2^{1 - \binom{t}{2}}$. Suppose now that $p$ is such that $p|V|=2n$ and observe that $p= n q^{-t +O(1)}$. If we choose a random subset of $V$ by picking each $v \in V$ independently with probability $p$, the expected number of monochromatic potential cliques in this subset is, for $n = 2^{t/2} q^{3t/8 + o(t)}$,
$$p^t 2^{1 - \binom{t}{2}} N_t \leq q^{- t^2 + o(t^2)} n^t 2^{-\frac{t^2}{2} + o(t^2)} q^{\frac{5t^2}{8} + o(t^2)} 
= \left(2^{-\frac{t}{2}} q^{-\frac{3t}{8} + o(t)} n\right)^t <1/2.$$
Since our random subset will also contain more than $n$ elements with probability at least $1/2$, there exists a choice of coloring and a choice of subset of order $n$ such that there is no monochromatic potential clique in this subset. This completes the proof. 

\vspace{3mm}
\noindent
{\bf Remark.}
Our method also gives a construction which matches Erd\H{o}s' bound $r(t) > \sqrt{2}^t$ up to lower-order terms. To see this, we set $V = \mathbb{F}_2^{2t}$ and color edges red or blue depending on whether $u \cdot v = 0$ or $1 \bmod{2}$. If we then sample $2^{t/2 + o(t)}$ vertices of $V$ at random, we can show that w.h.p.~the resulting set does not contain a monochromatic clique of order $t$. We believed this to be new, but, after the first version of this article was made public, we learned that such a construction was already discovered by Pudl\'ak, R\"odl and Savick\'y~\cite{PRS88} in 1988.
It was also pointed out to us by Jacob Fox that one can achieve the same end by starting with any pseudorandom graph on $n$ vertices for which the count of cliques and independent sets of order $2c \log_2 n$ is approximately the same as in $G(n, 1/2)$ and sampling $n^c$ vertices. This can be applied, for instance, with the Paley graph. 

\vspace{3mm}
\noindent
{\bf Acknowledgements.} We are extremely grateful to Vishesh Jain and Wojciech Samotij for reading an early draft of this paper and offering several suggestions which improved the presentation. We also owe a debt to Noga Alon and Anurag Bishnoi, both of whom pointed out the constraint on the dimension of isotropic subspaces, thereby improving the bound in our original posting.


\begin{thebibliography}{}

\bibitem{AK97}
N. Alon and M. Krivelevich, Constructive bounds for a Ramsey-type problem, {\it Graphs Combin.} {\bf 13} (1997), 217--225.

\bibitem{Con09}
{D. Conlon,} A new upper bound for diagonal Ramsey numbers, {\it Ann. of Math.} {\bf 170} (2009), 941--960.

\bibitem{Erd47}
P. Erd\H{o}s, Some remarks on the theory of graphs, {\it Bull. Amer. Math. Soc.} {\bf 53} (1947), 292--294. 

\bibitem{ESz35}
P. Erd\H{o}s and G. Szekeres, {A combinatorial problem in geometry,} {\it Compos. Math.} {\bf 2} (1935), 463--470.

\bibitem{L87}
H. Lefmann, A note on Ramsey numbers, {\it Studia Sci. Math. Hungar.} {\bf 22} (1987), 445--446.


\bibitem{PRS88}
P. Pudl\'ak, V. R\"odl and P. Savick\'y, Graph complexity, {\it Acta Inform.} {\bf 25} (1988), 515--535.

\bibitem{Sah20}
A. Sah, Diagonal Ramsey via effective quasirandomness, preprint available at arXiv:2005.09251 [math.CO].

\bibitem{Spe75}
J. Spencer, Ramsey's theorem --- a new lower bound, {\it J. Combin. Theory Ser. A} {\bf 18} (1975), 108--115.


\end{thebibliography}
\end{document}